\documentclass[11pt]{amsart}

\usepackage{mathpazo}
\usepackage[scaled=.95]{helvet}
\usepackage{courier}
  \usepackage{amsmath}
\usepackage{color}
\advance \topmargin by -\headheight
\advance \topmargin by -\headsep
\evensidemargin \oddsidemargin
\theoremstyle{plain}
\newtheorem{thm}{Theorem}[section]

\newtheorem{question}{Question}
\newtheorem{problem}{Problem}

\theoremstyle{definition}

\theoremstyle{remark}

\newcommand{\bR}{{\mathbb{R}}}

\newcommand{\bC}{{\mathbb{C}}}
\title{A Survey of Group Actions on $4$-Manifolds}
\author{Allan L. Edmonds}
\address{Department of Mathematics, Indiana University, Bloomington, IN 47405}
\email{edmonds@indiana.edu}
\begin{document}
\begin{abstract}
Over 50 years of work on group actions on $4$-manifolds, from the 1960's to the present, from knotted fixed point sets to Seiberg-Witten invariants, is surveyed. Locally linear actions are emphasized, but differentiable and purely topological actions are also discussed. The presentation is organized around some of the fundamental general questions that have driven the subject of compact transformation groups over the years and their interpretations in the case of $4$-manifolds. Many open problems are formulated. Updates to previous problems sets are given. A substantial bibliography is included.
\end{abstract}
\maketitle
 \setcounter{tocdepth}{1}
\tableofcontents
\section{Introduction}
These notes survey the realm of topological, locally linear actions of finite or compact Lie groups on topological 4-manifolds. 

The world of four-dimensional topological manifolds lies at the interface between the algebra of general $n$-manifolds, $n\ge 5$, and the geometry of $2$- and $3$-manifolds. Within that world of $4$-dimensions we now see a gradation from geometric and especially symplectic manifolds on the geometric side  and topological manifolds on the algebraic side, with pure differentiable (or, equivalently, piecewise linear) manifolds residing at the interface. Naturally this pattern persists in the study of group actions on $4$-manifolds.

The emphasis here is on highlighting those aspects of the higher dimensional framework where things go somewhat differently in dimension four and on highlighting positive topological results that pose challenges for differentiable transformation groups. Of special interest will be examples of topological locally linear actions on smooth manifolds that are not equivalent to smooth actions. We also indicate a few open problems related to possibly interesting phenomena related to purely topological actions.

For a particular group $G$ it is natural to study a hierarchy of group actions:  Free actions, which generalize both to semifree actions and to pseudofree actions, and then arbitrary actions.  In this survey all group actions are understood to be effective, in the sense each nontrivial group element acts nontrivially. We will generally consider only locally linear actions, except where explicitly stated otherwise. But, as we shall see, some of the  intriguing problems involve questions of whether particularly strange non-locally linear actions exist. We will also only consider orientation preserving actions unless explicitly stated otherwise, most notably in the examination of actions on the $4$-sphere.

The story of actions by compact connected Lie groups was largely developed prior to the explosion of 4-manifold topology in the 1980's. It is now a good time to look back and identify problems left over from that era. In particular there are issues about circle actions on non-simply connected $4$-manifolds that still deserve attention. Similarly some of the early work on the $4$-sphere, including knotted $2$-sphere fixed point sets and also free actions occurred early on and merits a new look.

The description of finite group actions is most successful in the world of simply connected topological 4-manifolds, which were completely classified 30 years ago in the work of Freedman \cite{FreedmanQuinn1990}.

Some of the more tantalizing problems, and problems on which there has been the most recent progress, include questions about smooth actions on smooth $4$-manifolds, addressed by methods of gauge theory, both Donaldson's Yang-Mills instantons and the more recent Seiberg-Witten theory of monopoles.

The results are more fragmentary and speculative in the world of non-simply connected 4-manifolds, including issues related to complements of $2$-dimensional fixed point sets, where most experts expect us eventually to find a failure of high dimensional surgery theory in the realm of 4-dimensional non-simply connected $4$-manifolds.

In its broad outlines this report  follows a historical progression, but we also attempt simultaneously to group topics together in a natural way, so that when appropriate we depart from the historical order of things.

The exposition is also organized around a sequence of general, often familiar, open-ended, guiding \textbf{questions}, many of which make sense in all dimensions, and a parallel sequence of  \textbf{problems}, which propose specific, concrete questions or conjectures, which are to the best of my knowledge still at least partially unsolved.

There are two appendices. In the first we list the problems from the classic Kirby Problem List (original 1976, updated 1995) related to group actions on $4$-manifolds, and update them to 2015. In the second appendix we list and update the problems on $4$--manifolds from the 1984 Transformation Groups Conference at the University of Colorado.

To accompany this paper we have prepared what we hope is a rather complete bibliography of papers on the subject of group actions on $4$-manifolds, including items not necessarily directly discussed in the body of the paper.

{\color{black}
The reader's attention is directed to a survey oriented toward group actions on symplectic 4-manifolds by W. Chen \cite{Chen2010}. Symplectic topology has seen substantial growth in recent years, so this survey is especially welcome. Although its content is largely distinct from the topics considered here, nevertheless we mention aspects of group actions on symplectic manifolds in several places.
}

\subsection*{Acknowledgements}
{This paper is an expanded and updated version of a lecture presented at the conference on Transformation Groups in Topology and Geometry, University of Massachusetts Amherst in 2008. The focus was on topological aspects of the subject. Thanks to participants in the conference for corrections, additions, and updates. {At the conference, M. Furuta gave a parallel survey on Gauge Theories and  Group Actions. I thank him for sharing the slides prepared for his lecture. The present brief overview of this aspect of the subject cannot compare with his more complete and informed presentation.} The text was further updated in subsequent years to incorporate comments and references provided by several people, as noted in the text, as well as to take into account later developments in the field. Thanks especially go to B.~Farb,  B.~Zimmermann, M.~Szymek, Q.~Khan, and, especially, to I.~Hambleton,  for their comments.} Note that some of the numbering of problems, in particular, has changed since the first version of the paper.

\section{Circle and Torus Actions}

Such actions of positive dimensional Lie groups are more accessible since the orbit space is a manifold of lower dimension. We focus on two questions:

\begin{question}[Equivariant Classification]
What is the classification of compact Lie group actions on $4$-manifolds, up to equivariant homemorphism, especially for the circle and torus groups?
\end{question}

This question has been satisfactorily answered through the combined work of Orlik and Raymond  \cite{OrlikRaymond1970,OrlikRaymond1974}, Fintushel  \cite{Fintushel1976,Fintushel1977,Fintushel1978}, Orlik \cite{Orlik1982}, Melvin and Parker  \cite{MelvinParker1986}.

The conclusion of all these cases is that one describes the action in terms of a weighted  orbit space, with the weight information describing the non-principal orbits and appropriate bundle data.  The details can be intricate.

\begin{question}[Topological Classification]
What $4$-manifolds admit actions by compact Lie groups of positive dimension, especially the torus or circle group, and how are they distinguished by their actions?
\end{question}

This question has pretty good answers in the case of simply connected $4$-manifolds or in the case of non-abelian  Lie groups. Results are still somewhat incomplete in the non-simply connected case for circle and even torus actions.

Melvin and Parker \cite{MelvinParker1986} have effectively handled the case of actions by non-abelian  Lie groups and reduced the general case to that of actions of the circle group $S^{1}$  and the torus $T^{2}=S^{1}\times S^{1}$. Building on earlier work of Melvin \cite{Melvin1981,Melvin1982},  Parker,  Orlik \cite{Orlik1982} and even earlier work of Richardson \cite{Richardson1961}, they prove that a $4$-manifold admitting an action of a non-abelian compact Lie group must be one of five types:  
\begin{enumerate}
\item $S^{4}$ or $\pm \mathbb{C}P^{2}$
\item connected sums of copies of $S^{1}\times S^{3}$ and $S^{1}\times \mathbb{R}P^{3}$
\item $SU(2)/H$-bundles over $S^{1}$, where $H$ is a finite subgroup of $SU(2)$
\item $S^{2}$-bundles over surfaces
\item certain quotients of $S^{2}$-bundles over surfaces by involutions.
\end{enumerate}

 Orlik and Raymond  \cite{OrlikRaymond1970,OrlikRaymond1974} had identified the simply connected 4-manifolds admitting torus actions, showing that they are connected sums of standard building blocks $S^{4}$, $\pm\mathbb{C}P^{2}$, and $S^{2}\times S^{2}$. Pao  \cite{Pao1977a,Pao1977b} extended this by showing that a general $4$-manifold with torus action is a connected sum of these building blocks, plus members of three other less-familiar  families in general (modulo the then undecided 3-dimensional Poincar\'e Conjecture). M.H. Kim  \cite{Kim1993} studied torus actions on non-orientable $4$--manifolds. 
 
Then Fintushel  \cite{Fintushel1978}  was able to identify completely the underlying manifolds admitting circle actions in  the simply connected case. Similar results are due independently  to Yoshida  \cite{Yoshida1978}. Fintushel used some of Pao's ideas in \cite{Pao1978} to show that the underlying manifold admits a possibly different circle action that extends to a torus action, must again be a connected sum of copies of $S^{4}$, $\pm\mathbb{C}P^{2}$, and $S^{2}\times S^{2}$ (again modulo the 3-dimensional Poincar\'e Conjecture). In particular this gives rise to many non-standard actions on standard manifolds, some of which we will allude to subsequently in the case of actions on $S^{4}$. Also the ``exotic'' topological $4$-manifolds, such as the $E_{8}$ manifold or the fake complex projective plane, etc., do not admit locally linear circle actions.

 Pao's replacement trick was to show how one can simplify the weighted orbit space (and the equivariant homeomorphism type of the action) while keeping the total space the same. As a first application this led to more knotted 2-spheres as fixed point sets in the $4$-sphere. The Pao trick has been applied more recently by Baldridge in a study of Seiberg-Witten invariants of 4-manifolds with circle action, which we'll mention again later in this survey.

\begin{problem}
Give a meaningful topological classification of the non-simply connected $4$-manifolds that admit $S^{1}$ actions.
\end{problem}

A special case of this was proposed by Melvin for the Kirby Problem List. See Appendix A. In particular, if the quotient manifold has free fundamental group, is the $4$--manifold a connected sum of copies of $S^{1}\times S^{3}$, $S^{2}\times S^{2}$, and $S^{2}{\otimes}S^{2}$ (the nontrivial bundle)?

{\color{black} Chen \cite{ChenArXiv2013} has studied the role of the fundamental group for $4$-manifolds with a fixed point free circle action. He shows among other things that for any finitely presented group with infinite center, there are at most finitely many distinct smooth (resp. topological) 4-manifolds which support a fixed-point free smooth (resp. locally linear) $S^1$-action and realize the given group as the fundamental group.

Smooth structures seem to play an interesting role in some of these issues about fundamental groups.  Chen \cite{ChenArXiv2011} shows that the homotopy class of an orbit of a smooth fixed point free circle action on smooth 4-manifold with a nonzero Seiberg-Witten invariant necessarily has infinite order in the fundamental group. Earlier Kottschick \cite{Kottschick2006} had shown that a closed symplectic $4$-manifold does not admit a smooth free circle action with contractible orbits using Seiberg-Witten theory.}

\begin{problem}
Determine how much of the theory of smooth or locally linear actions of the circle or the torus depends only on algebraic topology. How much of the analysis carries over to topological circle or torus actions or actions on homology manifolds or on Poincar\'e duality spaces, for example?
\end{problem}

Work of Huck \cite{Huck1995} and Huck and Puppe  \cite{HuckPuppe1998} is in this direction. The main goal is to use algebraic topology to prove that the intersection form of a 4-manifold with arbitrary circle action, defined on $H^{2}/\text{torsion}$ must split as an orthogonal sum of one and two-dimensional forms. This has been accomplished, generalizing the work of Fintushel in the simply connected, locally linear, case, but only provided the action satisfies the following property, 
which is a weakening of local linearity: Any fixed point has an invariant neighborhood containing at most four distinct orbit types.
This  remaining issue is  the main impediment to a satisfactory study of circle actions in the non-locally linear case.

\begin{problem}
Does a purely topological circle action on a $4$-manifold have at most four distinct orbit types in a small invariant neighborhood of any fixed point?
\end{problem}

Kwasik and Schultz identified this basic problem and achieved partial results, but the final result seems elusive now. One should note that the coned-off $E_{8}$-plumbing manifold admits a circle action that does have more than four orbit types in a neighborhood of the cone point. Does the closed $E_{8}$ topological $4$--manifold provided by Freedman's classification of simply connected topological $4$--manifolds admit a non-locally linear circle action?

Closely related to these considerations are the following two problems of Schultz from the Carlsson \cite{Carlsson1990} problem list of 1990.

\begin{problem}
Let $M$ be a closed topological 4-manifold with a topological circle action. 
Is the Kirby-Siebenmann triangulation obstruction of $M$ trivial? 
\end{problem}

\begin{problem}
Is every topological circle action on a 4-manifold concordant to a smooth 
action? 
\end{problem}
More recent work in the area of smooth circle actions has been addressed toward Seiberg-Witten invariants, which were calculated by Baldridge (\cite{Baldridge2001,Baldridge2003,Baldridge2004}). No applications to the underlying questions above were addressed, however. 

Herrera  \cite{Herrera2006} has shown that if a smooth $4$-manifold  has an even intersection pairing, but is not necessarily spin, and admits a nontrivial smooth circle action, then the signature necessarily vanishes. This extends in dimension four the famous result of Atiyah-Hirzebruch that the $\widehat{A}$ genus vanishes for smooth spin manifolds admitting a smooth circle action. In brief Herrerra shows that the circle action lifts to a circle action on a spin covering.

Finally, there is a whole additional thread of symplectic circle and torus actions in dimension four that we have not been able to address here, for reasons of lack of time and expertise.

\section{The $4$-sphere}
Changing focus to finite group actions, a natural place to begin is with actions on the $4$-sphere.

\subsection{Codimension one fixed set}
Poenaru  \cite{Poenaru1960} and, independently, Mazur   \cite{Mazur1961} were the first to observe non-simply connected homology spheres that bound contractible manifolds whose doubles are spheres, thus producing the first non-linear actions\footnote{Except for the suspension of Bing's (1952) involution on $S^{3}$ fixing a wildly embedded $2$-sphere.} on $S^{4}$. DeRham    \cite{deRham1965} elaborated Poenaru's construction to produce infinitely many examples.

Given Freedman's work on topological $4$-manifolds \cite{FreedmanQuinn1990} we have a complete understanding in the topological category. Such topological, locally linear, reflections (clearly orientation-reversing) are in one-to-one correspondence with their fixed point sets, which are integral homology $3$-spheres  bounding contractible topological $4$-manifolds.

Applying this construction to a homology $3$-sphere that does not bound a smooth contractible manifold, one obtains non-smoothable but locally linear actions on $S^{4}$.

\subsection{Smith conjecture}
The next earliest results refer to 2-dimensional fixed point sets.
\begin{question}[Generalized Smith Conjecture]
Can a knotted $2$-sphere be the fixed point set of a $C_{p}$ action on $S^{4}$.
\end{question}
Giffen  \cite{Giffen1966} gave the first examples of such knotted fixed point sets, via cyclic branched covers of certain twist spun knots of Zeeman, valid in all dimensions greater than $3$.  The first examples were restricted to odd order groups.
Gordon  \cite{Gordon1974} extended the construction to cover any order cyclic group. Then Pao  \cite{Pao1978} used his ``replacement trick'' to give a complete analysis (modulo the classical Poincar\'e Conjecture) of the situation for circle actions, and in particular gave infinite famiies of examples.

\begin{problem}
Identify more obstructions to a knot in $S^{4}$ being a fixed point set.
\end{problem}
What if anything meaningful can one say about the classical invariants of fixed point knots in this dimension?
It is a theorem of Kervaire, valid in the smooth, PL, or topological locally flat categories, that
all $2$-knots in $S^{4}$ are slice.  McCooey \cite{McCooey2007c} has shown that a  $2$-sphere fixed by a  $C_{p}$ action is equivariantly slice, in the topological locally linear category

\subsection{Free actions} According to the Lefschetz Fixed Point Theorem the only group acting freely on a $4$-dimensional sphere is $C_{2}$. Moreover, the orbit space of such an action is homotopy equivalent to $\mathbb{R}P^{4}$. 

\begin{question}[Fixed Point Free Involutions]
Is there a fixed point free involution on $S^{4}$ not equivalent to the antipodal map?
\end{question}

Cappell and Shaneson  \cite{CappellShaneson1976}  produced the first examples of smooth $4$-manifolds homotopy equivalent but not diffeomorphic to $\mathbb{R}P^{4}$, by an interesting construction removing a neighborhood of a torus and sewing back in something else. For some while the universal coverings of these examples held promise of providing counterexamples to the smooth four-dimensional Poincar\'e conjecture.  Eventually, however, Gompf  \cite{Gompf1991} showed that one of these examples has universal covering diffeomorphic to $S^{4}$. 
{\color{black}
More recently Akbulut \cite{Akbulut2010} has shown that all the candidates in an infinite list of Cappell and Shaneson spheres are diffeomorphic to one another. 
It follows that none of these manifolds is a counterexample to the smooth four-dimensional Poincar\'e conjecture. And Gompf \cite{Gompf2010} has further shown that the members of an even larger family of examples are standard. 
}

Meanwhile Fintushel and Stern  \cite{FintushelStern1981} had used different techniques to produce the first smooth fake $\mathbb{R}P^{4}$ with universal covering diffeomorphic to $S^{4}$.

Ruberman  \cite{Ruberman1984} showed how to construct the unique topological  4-manifold homotopy equivalent to $\mathbb{R}P^{4}$ but not homeomorphic to it, as promised by surgery theory.

\begin{problem}[Ruberman]
Does any smooth, free involution on $S^{4}$ admit an invariant $2$-sphere?
\end{problem}
Of course, the standard antipodal map has an invariant $2$-sphere, and also Ruberman's topological construction provides just such an invariant (topological) $2$-sphere. At least some of the constructions mentioned above do yield an invariant $2$-sphere.

These studies of free involutions predate the introduction of gauge theory into 4-dimensional topology. Are there useful contributions to these questions that might arise from gauge theory?
\subsection{One fixed point actions}
\begin{question}
Can a finite group act locally linearly on $S^{4}$ with just one fixed point?
\end{question}
Without the explicit local linearity hypothesis, this would be equivalent, via one-point compactification, to actions on $\mathbb{R}^{4}$ without a fixed point, which we will also consider separately below. Of course, by Smith Theory, any $p$-group fixes at least two points, so the full singular set for the action consists of more than just a single fixed point.

In an interesting turn of events, this was originally answered ``No''  in the smooth case via gauge theory, by Furuta  \cite{Furuta1989}, in the first application of Donaldson's theory of (anti) self-dual connections to transformation groups. Only later was a topological, locally linear version proved, by DeMichelis  \cite{DeMichelis1989}, using more traditional techniques. What about a purely topological action? One can argue, using Smith theory, that the answer is still ``no'' for arbitrary topological actions of solvable groups. (In higher dimensions this argument breaks down and there do exist smooth periodic maps of non-prime order on spheres with just one global fixed point.)
\begin{problem}
Can a finite, non-solvable group act topologically on $S^{4}$ with just one fixed point?
\end{problem}
We remark that there is indeed a one-fixed-point action of the group $A_{5}$ on a $4$-dimensional integral homology manifold of the homotopy  type of $S^{4}$. For there is a smooth action on the Poincar\'e homology  $3$-sphere with just one fixed point, leading to a one-fixed-point action on its reduced suspension.
\subsection{Pseudofree actions} For our purposes an action of a finite group on a manifold is said to be \emph{pseudofree} if it is free away from a discrete set of singular points with nontrivial isotropy groups.

Kulkarni  \cite{Kulkarni1982} studied pseudofree actions of general finite groups on general spaces. He largely classified the finite groups that can act pseudofreely, preserving orientation, on a space with the homology of an even dimensional homology sphere. He left open whether a dihedral group $D_{n}$ of order $2n$, $n$ odd, could act pseudofreely on an even dimensional sphere of dimension greater than $2$.

This author \cite{Edmonds2010} subsequently completed the argument that the dihedral group $D_{n}$ does not act pseudofreely, preserving orientation, on $S^{4}$, or, more generally, on $S^{d}$, where the dimension $d$ is divisible by 4. Then  Hambleton \cite{Hambleton2011} was able to apply rather different techniques to show further that $D_{n}$ does not act pseudofreely, preserving orientation, on any even-dimensional sphere.

The study of group actions with isolated fixed points has also generated some work related to topological actions that are not necessarily nice near the singular points. Kwasik and Schultz  \cite{KwasikSchultz1990} studied orientation-preserving, pseudofree actions of cyclic groups on $S^{4}$ in both the locally linear and topological categories, exploring the projective class group obstruction to desuspending an action of $C_{2n}$ on $S^{4}$ to be the join of a free action on $S^{3}$ with a nontrivial action on $S^{0}$.
See also the general discussion of conelike and weakly conelike  actions in Freedman-Quinn \cite{FreedmanQuinn1990}, Section 11.10. They argue, in particular, that the weakly conelike actions are precisely the pseudofree actions whose associated projective class group obstruction vanishes.

\subsection{Groups acting on $S^{4}$}
Here we consider group actions with no particular assumptions about the fixed point sets.
\begin{question}
What finite groups act on $S^{4}$?
\end{question}

Of course, if the group is locally linear and has a fixed point, then it is a subgroup of $O(4)$, and if it acts linearly, then it is a subgroup of $O(5)$.

\begin{problem}
Show that a finite group acting on $S^{4}$, without fixed point, is isomorphic to a subgroup of $O(5)$.
\end{problem}
{\color{black}
Mecchia and Zimmermann  \cite{MecchiaZimmermann2011} came close to solving this problem in the case of smooth actions, building on their earlier work in \cite{MecchiaZimmermann2006}, which addressed the case of non-solvable groups. The orientation-preserving case has finally been completed by Chen, Kwasik, and Schultz \cite{ChenKwasikSchultzArXiv2014} dealing with a delicate issue of subgroups of index 2. The orientation-reversing case is not completely settled. These results should apply equally well to both smooth and locally linear actions. What if the local linearity is dropped? Then Chen, Kwasik, and Schultz \cite{ChenKwasikSchultzArXiv2014} construct a \emph{topological} orientation-reversing action of a group that does not embed in $O(5)$.
}

\subsection{The Smith conjecture on representations at fixed points}
\begin{question}
If a group $G$ acts locally linearly (or smoothly) on a sphere with two fixed points, are the local representations at the two fixed points are equivalent?
\end{question}
In higher dimensions the answer is no in general and the problem has provided a rich line of research. In dimension $4$, in the smooth case, an affirmative answer was given by Hambleton and Lee   \cite{HambletonLee1992}  and Braam and Matic \cite{BraamMatic1993} in another one of the early applications of gauge theory to group actions. 

\begin{thm}[Hambleton and Lee   \cite{HambletonLee1992}, Braam and Matic \cite{BraamMatic1993}]
If a finite group acts smoothly on $S^{4}$ with two isolated fixed points, then the local tangential  representations at the fixed points are equivalent.
\end{thm}

The idea is to look at the instanton number one moduli space of self-dual connections with its induced action and relate the normal representations to the fixed points on the original manifold to the normal representations to fixed points in the moduli space near a reducible connection. The details are subtle.
On the one hand, their proofs illustrated the fact that a gauge-theoretic approach can be intuitively appealing and eliminate a study of numerous special cases.  On the other hand the gauge-theoretic arguments have their own subtleties, related to the failure of equivariant transversality, and were hard to get right. See Section \ref{sec:gauge} for a little more discussion.

DeMichelis  \cite{DeMichelis1989} was able to give a proof valid for topological, locally linear actions on homology $4$-spheres, using more classical techniques, after learning of the early attempts to prove this result via gauge theory.
\section{Euclidean $4$-space}
In some ways this might have been the appropriate place to start our survey of finite group actions, but non-compactness raises its own difficulties.
\begin{question}
Is there a periodic self-map of $\mathbb{R}^{4}$ that has no fixed point?
\end{question}
In higher dimensions such actions were first found by Conner and Floyd, with examples found in all dimensions $\ge 8$ by Kister. Smith  \cite{Smith1960} observed that no such periodic maps exist in dimensions less than 5 (and not in dimension 5 or 6 either if one assumes the maps are smooth). Examples were found in dimension $7$ by Haynes, et al. \cite{HaynesKwasikMastSchultz2002}. Smith's argument extends easily to the case of solvable groups. Thus it remains to contemplate actions of more complicated groups.

\begin{problem}
Find a non-solvable group acting without a global fixed point on $\mathbb{R}^{4}$.
\end{problem}

Contemplating possible actions by unfamiliar groups on $\mathbf{R}^{4}$, one is led to ask just what groups arise.
\begin{problem}
Show that a finite group that acts on $\mathbb{R}^{4}$ is isomorphic to a subgroup of $O(4)$.
\end{problem}
Of course, if the action is smooth or locally linear and has a fixed point, then the group is a subgroup of $O(4)$. 
{\color{black}This problem has been completely resolved in the smooth or locally linear cases by Guazzi, Mecchia, and Zimmermann
 \cite{GuazziMecchiaZimmermann2011}, without, however, proving that such a group action necessarily has a fixed point. Indeed their argument applies to actions on acyclic manifolds and there is an action of the alternating group $A_{5}$ on such an acyclic $4$-manifold with no global fixed point. On the other hand a surgery-theoretic analysis of Hambleton and Madsen \cite{HambletonMadsen1986}, appropriately translated into dimension $4$, shows that there exist groups $G$ that act semifreely on $\mathbb{R}^{4}$ with a single isolated fixed point such that $G$ is not isomorphic to a subgroup of $O(4)$. See the discussion in Hambleton's recent survey \cite{Hambleton2015} after the statement of Theorem 9.1 for a brief, but explicit, discussion of extending surgery theory finite groups to dimension $4$. He also mentions some other sorts of exotic actions that arise from this point of view, including examples of non-solvable groups that act topologically, but not linearly.}
 
Musing about the ``fake'' $\mathbb{R}^{4}$'s that came out of the work of Freedman and Donaldson in the 1980's one is led to wonder about their group actions.

\begin{problem}
Is there a smooth fake $\mathbb{R}^{4}$ that admits no nontrivial, smooth group  actions?
\end{problem}
This problem was already raised at the 1984 Transformation Groups conference in Boulder. Taylor \cite{Taylor1999} has constructed examples of fake  $\mathbb{R}^{4}$s on which the circle group cannot act smoothly.

\section{Simply connected $4$-manifolds}
As the study of group actions on $4$-manifolds developed, interest began to move beyond the basic spheres, disks, and euclidean spaces associated with Smith theory, to look more closely at standard $4$-manifolds including the complex projective plane $\mathbb{C}P^{2}$, $S^{2}\times S^{2}$, and the K3 surfaces of algebraic geometry. Along the way a rather successful study of existence and classification of locally linear, pseudofree actions on closed, simply connected $4$--manifolds ensued.

\subsection{Particular model $4$-manifolds}
After the sphere, perhaps the most important basic closed $4$-manifold is the complex projective plane. It certain ways it  is even  more important than the sphere.
\subsubsection{The complex projective plane}
\begin{question}
What groups act on $\mathbb{C}P^{2}$?
\end{question}
There is an excellent answer.
\begin{thm}[Wilczy\'nski \cite{Wilczynski1987,Wilczynski1990}, Hambleton and Lee \cite{HambletonLee1988}]
If a finite or compact Lie group acts locally linearly  on $\mathbb{C}P^{2}$, then $G$ admits a faithful projective representation of complex dimension $3$, and hence also acts linearly on $\mathbb{C}P^{2}$.
\end{thm}
The first version, for actions of finite groups that are homologically trivial was due independently to Hambleton and Lee \cite{HambletonLee1988} and Wilczy\'nski\cite{Wilczynski1987}. In that case one can just say that $G$ must be a subgroup of $PGL_{3}(\mathbb{C})$. The general result is due to Wilczy\'nski \cite{Wilczynski1990}.  He characterizes the ``universal group'' acting as the quotient $NU(3)/ZU(3)$, where $NU(3)$ and $ZU(3)$ denote the normalizer and the centralizer, respectively, of the unitary group $U(3)$ when viewed as a subgroup of $O(6)$.

Wilczy\'nski \cite{Wilczynski1988} also gave a similar result for groups acting locally linearly on the Chern manifold, the unique $4$-manifold homotopy  equivalent to $\mathbb{C}P^{2}$ but not homeomorphic to it. In that case the possible groups, however, are more restrictive, since $G$ is necessarily finite and without $2$-torsion. In this case it also follows that any action is necessarily pseudofree. For more general $4$-manifolds with the homology of $\mathbb{C}P^{2}$, any group of symmetries is known to be a linear group. 

\begin{question}
Is every locally linear action on $\mathbb{C}P^{2}$ equivalent to a (projective) linear action?
\end{question}
There is just one known exception: Hambleton and Lee observed an action of $C_{p}\times C_{p}$, where the $2$-sphere fixed by one $C_{p}$ is knotted with respect to the fixed point set of another $C_{p}$.

For simplicity, especially in light of the preceding example, we restrict to semifree actions.  

\begin{thm}[Wilczy\'nski  \cite{Wilczynski1991}]
A locally linear, semifree action of a finite cyclic group on $\mathbb{C}P^{2}$ is equivalent to a linear action.
\end{thm}

According to Edmonds and Ewing  \cite{EdmondsEwing1989}, in the pseudofree case, the local representations at the three isolated fixed points must be the same as those of a linear action. (Dovermann had already shown that in the case the fixed point set consists of a $2$--sphere and an isolated point, then the local representation at the isolated point is standard.) The proof involved showing that the only solutions of the equations of the $g$-signature theorem are those coming from standard linear actions. Hambleton and Lee \cite{HambletonLee1992} gave an independent proof of this, assuming the action is smooth, as an application of equivariant gauge theory. Then a four-dimensional surgery argument shows that the given action is equivalent to the corresponding linear model. See Hambleton and Lee-Madsen  \cite{HambletonLeeMadsen1989} and especially Wilczy\'nski  \cite{Wilczynski1991} for full details of a more general result.

\subsection{General simply connected 4-manifolds}

According to Freedman, closed simply connected $4$-manifolds are classified up to homeomorphism by the intersection form and, in the case of intersection forms of odd type, the stable triangulation obstruction.

\subsubsection{Existence of group actions}\label{subsec:existence}
\begin{question}
What simply connected $4$-manifolds admit actions of $C_{p}$?
\end{question}

\begin{thm}[Kwasik and Vogel \cite{KwasikVogel1986b}]
If $M$ is a closed, orientable, $4$-manifold and $C_{2}$ acts locally linearly on $M$, then the Kirby-Siebenmann triangulation obstruction $\text{KS}(M)\in H^{4}(M;\mathbb{Z}_{2})$ must vanish.
\end{thm}

This extended a slightly earlier result of Kwasik \cite{Kwasik1986a}, which applied specifically to the Chern manifold, or fake $\mathbb{C}P^{2}$. At the same time Kwasik argued that the Chern manifold admits locally linear, pseudofree, actions of odd order cyclic groups.

\begin{thm}[Edmonds \cite{Edmonds1987}]
If $M$ is a closed, orientable, $4$-manifold and $p$ is a prime greater than $3$, then $C_{p}$ acts locally linearly, homologically trivially, and pseudofreely on $M$.
\end{thm}
These manifolds all admit locally linear actions of $C_{3}$, too, but in some cases a two-dimensional fixed point set is required.
\begin{question}
What automorphisms of the intersection form of a simply connected $4$-manifold $M$  are induced by actions of $C_{p}$ on $M$?
\end{question}

In the case of pseudofree actions Edmonds and Ewing gave a good answer for groups of prime order.

\begin{thm}[Edmonds and Ewing \cite{EdmondsEwing1992}] Let a closed, simply connected $4$--manifold $M$ be given, together with a representation of $C_{p}$ on $H_{2}(M)$, preserving the intersection form. Then there is a locally linear, pseudofree action of $C_{p}$ on $M$ inducing the given action on homology if and only if the representation on $H_{2}(M)$ is of the form $\mathbb{Z}[C_{p}]^{r}\oplus \mathbb{Z}^{t}$ and there is candidate fixed point data for $t+2$ isolated fixed points satisfying the $g$--signature theorem and an additional Reidemeister torsion condition.
\end{thm}
Note that the number of fixed points is determined by the action on homology by the Lefschetz fixed point theorem. 

This result has been useful for showing that certain kinds of actions exist topologically, even as it has been shown that the action cannot be smooth.

The construction of such actions proceeds by building up a smooth equivariant handlebody structure, with one $0$--handle and suitable $2$--handles. One then faces the problem of capping off the manifold with boundary, which is necessarily an integral homology $3$--sphere with a free action of $C_{p}$. It is here that the signature and torsion condition comes in as one applies high-dimenionsal surgery theory ideas to a $4$--dimensional situation with finite fundamental groups to show that the action on the boundary is equivariantly $h$--cobordant to a linear action on $S^{3}$. And it is here that  one is forced to use a topological as opposed to a smooth or PL manifold.

{\color{black}
\begin{problem}
Develop an analog of the Edmonds and Ewing existence theorem \cite{EdmondsEwing1992} for actions that also allows for fixed point sets including $2$-spheres or, even, other surfaces of positive genus, not just isolated points.
\end{problem}
}

{\color{black}
Even in the case of homologically trivial actions this may not be so easy. If the group $C_{p}$ acts locally linearly and homologically trivially on a closed, oriented, simply connected $4$-manifold $M$, then an  Euler characteristic argument shows that the number $n$ of $2$-spheres in the fixed point set satisfies $n\le \chi(M)/2$. One may ask what values of $n$ are realizable. Recently Hamilton \cite{Hamilton2016} has used a close analysis of aspects of the $G$-Signature Formula to show that in particular cases there are more subtle restrictions. For example if all $2$-spheres in the fixed point set have self-intersection $\le s<0$, where $s$ is a negative integer, then
\[
n < \frac{\chi(M)}{2-s}\left(1 + \frac{6}{p(2-s)-(4+s)}\right).
\]
If the the manifold $M$ is smooth and minimal (no $(-1)$-spheres) and has $b_{2}^{+}>0$ and a nonzero Seiberg-Witten invariant (implying that self-intersections of smooth $2$-spheres are negative) then this implies
\[
n < \frac{\chi(M)}{4}\left(1 + \frac{3}{2p-1} \right).
\]
}
{\color{black}
\subsubsection{Groups that act on some $4$-manifold}
Of course \emph{any} finite group is the fundamental group of a closed, orientable, smooth $4$-manifold and hence acts freely on the universal covering of that manifold. But if the group acts homologically trivially on a manifold more complicated than the $4$-sphere, $\mathbb{C}P^{2}$ or $S^{2}\times S^{2}$, then there are significant restrictions on the finite groups that occur.

A first result in this direction is due to Hambleton and Lee \cite{HambletonLee1995} who showed, making use of gauge theory, that a finite group that acts smoothly on a connected sum of $n\ge 1$ copies of the complex projective plane,  inducing the identity on homology, is isomorphic to a subgroup of $PGL_{3}(\mathbb{C})$, and that if $n>1$, then in fact the group is abelian of rank $2$. They also analyzed the local representations at fixed points and showed as a consequence of work of Edmonds and Ewing \cite{EdmondsEwing1992} in the topological category that there are smooth cyclic group actions on a connected sum of $n\ge 1$ copies of the complex projective plane,  inducing the identity on homology, that are not smoothable.

More generally McCooey gave the following general result in the topological locally linear context.

\begin{thm}[McCooey \cite{McCooey2002a}]
Let $G$ be a compact Lie group acting effectively, locally linearly, and homologically trivially on a closed $4$-manifold $M$, with $H_1(M;\bold Z)=0$. If  $b_2(M)\geq 3$, then $G$ is isomorphic to a subgroup of $S^1 \times S^1$. 
\end{thm}

McCooey shows that the same conclusion applies when $b_2(M)= 2$ and ${\rm Fix}(M,G) \neq \emptyset$. When $b_{2}(M)=2$, Mecchia and Zimmerman \cite{MecchiaZimmermann2009} prove in general that the only simple group that occurs is $A_5$, and also give a short list of possible finite nonsolvable groups that are candidates for actions of such groups. 

McCooey \cite{McCooey2007b} has also studied some cases where $H_1(M;\mathbb{Z})\neq 0$ and is torsion-free. Assuming in addition that $b_{2}(M)\neq0, 2$ and that $\chi(M)\neq 0$, then he shows that an orientation-preserving finite group acting on $M$ in a homologically trivial way must be cyclic.

Recently Mundet i Riera has a provided an alternative perspective in analogy with the classical Jordan theorem about orders of finite subgroups of $\text{GL}(n,\mathbb{R})$, as follows.
\begin{thm}[Mundet i Riera \cite{MundetiRieraArXiv2015}]
Let $M$ be a compact, orientable, connected $4$-dimensional
smooth manifold, 
such that $\chi(M)\neq 0$. There exists a natural number $C$ such that
any finite group $G$ acting smoothly and effectively on $M$ has an abelian
subgroup $A$ of index $\leq C$ with $\chi(M^A)=\chi(M)$. If $\chi(M)>0$ then $A$ has rank at most $2$, and if $\chi(M)<0$ then $A$ is cyclic.
\end{thm}
}

\subsubsection{Classification of group actions}
in the case of actions on standard manifolds one attempts to show that an arbitrary action is equivalent to a standard action. But for general manifolds one does not necessarily have a standard model action.
\begin{question}
What is the classification of actions of  a finite cyclic group $C_{n}$ on  a given simply connected $4$-manifold?
\end{question}
Hambleton and Kreck \cite{HambletonKreck1988,HambletonKreck1990,HambletonKreck1993c} addressed the case of free actions, by classifying the orbit spaces, culminating in the following result.

\begin{thm}[Hambleton and Kreck \cite{HambletonKreck1993c}] 
A closed, oriented $4$-manifold with finite cyclic fundamental group is classified up to homeomorphism by the fundamental group, the intersection form on $H_{2}(M;\mathbb{Z})/\text{Tors}$, the $w_{2}$-type, and the Kirby-Siebenmann triangulation obstruction. Moreover, any isometry of the intersection form can be realized by a homeomorphism.
\end{thm}
Actually Hambleton and Kreck only assume that the cohomology of $\pi_{1}$ is $4$--periodic, a condition that Bauer \cite{Bauer1988} weakened to requiring that the $2$--Sylow subgroup have $4$--periodic cohomology. For more general fundamental groups there are additional invariants, including a class in $\mathbb{Z}/|\pi_{1}|\mathbb{Z}$ and a secondary obstruction in $\text{Torsion}(\Gamma(\pi_{2}(X))\otimes \mathbb{Z})$, where $\Gamma$ denotes Whitehead's quadratic functor.

Wilczy\'nski  \cite{Wilczynski1994} then dealt with general pseudofree actions of cyclic  groups, completed by  Wilczy\'nski and Bauer \cite{BauerWilczynski1996}

Up to a finite ambiguity pseudofree actions on simply connected $4$-manifolds are determined by  the local representations at fixed points and by the action on homology. The most all-inclusive theorem would be the following.

\begin{thm}[Wilczy\'nski \cite{Wilczynski1994}, Wilczy\'nski and Bauer \cite{BauerWilczynski1996}
]
Suppose $M$ is a closed, oriented, simply-connected, topological $4$-manifold with an action of $G=C_{n}$, acting locally linearly and pseudofreely, preserving orientation, and with non\-empty  fixed point set $M^{G}$. If $n$ is odd or $w_{2}(M)\ne 0$, then the oriented topological conjugacy classes of locally linear, pseudofree $G$-actions on $M$, with the same local  fixed point representations and the same equivariant intersection form (and the same orbit space triangulation obstruction) are in one-to-one correspondence with the elements of a certain finite double coset space $O(\nu)\backslash O_{*}(\lambda)/O(\lambda)$, and are distinguished by the quadratic $2$-type.
\end{thm}
\subsection{Extension}According to Freedman's theory of topological $4$-manifolds \cite{FreedmanQuinn1990}, any homology $3$-sphere bounds a unique contractible topological $4$-manifold. 

\begin{problem}
Does every finite group action on a homology $3$-sphere extend to one on the contractible topological $4$-manifold it bounds?
\end{problem}
The answer is known to be yes in the case of free actions on homology $3$-spheres, provided one does not insist on local linearity at the central fixed point. Versions of this are due to Kwasik and Vogel \cite{KwasikVogel1986b}, Ruberman, and Edmonds \cite{Edmonds1987}. The problem of making a locally linear extension was dealt with by Edmonds \cite{Edmonds1987} and Edmonds and Ewing \cite{EdmondsEwing1992} as the last step in the construction of  locally  linear pseudofree actions on closed $4$--manifolds. See also Kwasik and Lawson \cite{KwasikLawson1993}
{\color{black}
who used gauge theory to show that certain free cyclic actions on Brieskorn spheres that bound smooth, contractible $4$-manifolds do not extend smoothly. Recently Anvari and Hambleton \cite{AnvariHambleton2014} have used gauge-theoretic considerations to show that none of the standard free actions on Brieskorn spheres extend smoothly over any contractible smooth manifold they may bound.
} 
An interesting remaining case is that of actions of $C_{p}$,  with nonempty fixed point set a  knotted circle, especially the case $p=2$. A related non-equivariant question would be whether a knot in a homology $3$-sphere bounds a disk in a contractible $4$-manifold, in analogy with the way that any knot in $S^{3}$ can be coned off in $D^{4}$. The analogous question for circle actions is closely related to the question mentioned earlier about whether a topological circle action can have more than four orbit types in a neighborhood of a fixed point.

\begin{problem}
Do all closed simply connected $4$-manifolds (even with nonvanishing  triangulation obstruction) admit topological involutions? Does the Chern manifold, the fake $\mathbb{C}P^{2}$, admit an involution?
\end{problem}

Such an involution cannot be locally linear when the triangulation obstruction is nonzero, by Kwasik and Vogel  \cite{KwasikVogel1986a,KwasikVogel1986b}. Kwasik and Vogel did show that some nontriangulable $4$--manifolds do admit involutions, including the case of closed, simply connected, spin $4$-manifolds.
\subsection{Other special manifolds}
\subsubsection{K3 surface} By a topological K3 surface we mean a closed, simply connected topological $4$-manifold with intersection pairing on second homology that is even of rank $22$ and signature $-16$.  By Freedman's classification all such manifolds are homeomorphic to a standard algebraic surface, such as the hypersurface in $\mathbb{C}P^{3}$ consisting of solutions to the equation $z_{0}^{4}+z_{1}^{4}+z_{2}^{4}+z_{3}^{4}=0$.
It follows that all such manifolds  are smoothable. From the point of view of algebraic geometry a K3 surface is a smooth, simply connected, complex algebraic variety with trivial canonical bundle. All such algebraic surfaces are diffeomorphic. But gauge theory has led to constructions of many \emph{homotopy K3 surfaces,} that is smooth manifolds homotopy equivalent, and hence homeomorphic, to the standard K3 surface, that are definitely not diffeomorphic to the standard K3 surfaces.

The many nontrivial algebraic K3 surfaces admit rich families of nontrivial symmetries, and there is a vast literature on the subject in algebraic geometry. We refer to Nikulin \cite{Nikulin1979}, who initiated much of the work in this area, to Mukai  \cite{Mukai1988} who classified the maximal groups acting symplectically (i.e., acting trivially on the one-complex-dimensional space of holomorphic $2$-forms), and for the most recent work in the area, including the non-symplectic case, to work of Sarti  and her coauthors, e.g.\cite{ArtebaniSarti2008, GarbagnatiSarti2007}.

There has been a variety of work aimed at identifying non-smoothable actions on such a surface.

\begin{question}
Do smooth homotopy $K3$--surfaces admit a smooth, homologically trivial group actions?
\end{question}
Complex analytic group actions are known to be homologically nontrivial. On the other hand the work of Edmonds  \cite{Edmonds1987} showed that such a manifold admits homologically trivial, locally linear actions. Thus the question is really about smooth actions. 
\begin{thm}[Ruberman \cite{Ruberman1995}, Matumoto \cite{Matumoto1992}]
A topological K3 surface does not admit a homologically trivial, locally linear, involution.
\end{thm}

According to Jin Hong Kim  \cite{Kim2008}  a homotopy K3 surface does not admit a homologically trivial smooth $C_{3}$ action.  The proof uses somewhat delicate computations involving Seiberg-Witten invariants. But see the comment in the appendix updating Problem 4.124 in the classic Kirby problem List.

\begin{problem}
Show that $C_{p}$, $p\ge 3$, cannot act smoothly and homologically  trivially on a homotopy $K3$ surface.
\end{problem}
 {\color{black}
Chen and Kwasik \cite{ChenKwasik2007a} have used Seiberg-Witten theory to show that a finite group acting homologically trivially on a \emph{symplectic} $4$-manifold with nonzero signature, $c_{1}^{2}= 0$, and $b_{2}^{+}\ge 2$ (e.g. a K3 surface) must act trivially.

Chen and Kwasik \cite{ChenKwasik2008} have also found examples of homotopy K3-surfaces that do not admit smooth actions of certain groups that do act on the standard K3-surface, by examining the action of the group permuting the Seiberg-Witten basic classes.
}

 \subsubsection{Other building block manifolds}ÊÊÊ
 
$E_{8}$: Edmonds \cite{Edmonds1997} analyzed the automorphisms of the $E8$ lattice and determined which ones could be realized by locally linear periodic maps.

$S^{2}\times S^{2}$: Wilczynski \cite{Wilczynski1994} showed that there exist inequivalent (but stably equivalent) actions that induce the same action on homology. McCooey \cite{McCooey2007a} has characterized the groups that can act pseudofreely, under the assumption of local linearity, as standard subgroups of the normalizer of $O(3)\times O(3)$ in $O(6)$.

$E_{4m}$: This is a family of definite intersection pairings for $4$-manifolds, that are spin if and only $m$ is even. Jastrzebowski  \cite{Jastrzebowski1995} determined which manifolds admit locally linear involutions with exactly two isolated fixed points, corresponding to the case where the second homology group is a free $\mathbb{Z}[C_{2}]$-module. When he combines his algebraic analysis with the results of Edmonds \cite{Edmonds1989}, one conclusion is that a closed, simply connected topological $4$-manifold with intersection pairing $E_{4m}$ admits a locally linear involution if and only if $m$ is odd and the Kirby-Siebenmann triangulation obstruction vanishes. In particular the $E_{16}$ $4$-manifold admits no locally linear involution, despite having a vanishing triangulation obstruction.

\subsection{Conjugation Spaces}
Hambleton and Hausmann \cite{HambletonHausmann2011} have studied these special actions of $C_{2}$ on connected 4-manifolds, achieving a very satisfactory classification. The definition is too technical to give here, but implies that the fixed point set is a connected surface and the orbit space is a  mod 2 homology 4-sphere. Indeed, their main theorem, shows that the fixed point set and orbit space together give a complete classification of these actions. The classical examples of such actions come from complex conjugation on $\bC P^{2}$ and $S^{2}\times S^{2}$.

\subsection{Permutation Representations}
By results of Edmonds \cite{Edmonds1989}, interpreting the Borel spectral sequence in equivariant cohomology in dimension $4$, non-permutation representations on the second homology of simply connected $4$--manifolds correspond to fixed surfaces of positive genus.
\begin{question}
For a $4$-manifold $M$, what representations on $H_{2}(M)$ are induced by group actions?
\end{question}
Of course the representation must preserve the intersection form.
\begin{problem}
Show  that an action of a finite group on a closed, simply connected 4-manifold, with definite intersection form, must induce a signed permutation representation on $H_{2}$.
\end{problem}

By Donaldson the intersection form of a smooth positive definite manifold is diagonalizable, and it follows that all automorphisms groups yield signed permutation representations automatically. Therefore this is only a question for topological $4$-manifolds. A detailed analysis of $\text{Aut} E_{8}$, where non-permutation representations  arise algebraically, showed that this problem has an affirmative answer in this one case. Every integral representation of $C_{2}$ is automatically a signed permutation representation. So the result is also a ``yes'' answer for this group.

In a related vein, work in Edmonds \cite{Edmonds2005} produces a $C_{25}$ action that is pseudofree but not semifree, while inducing a permutation representation.  But work of Hambleton and Tanase \cite{HambletonTanase2004} has shown that a \emph{smooth} pseudofree action, on a positive definite $4$--manifold, is semifree, via ASD Yang-Mills gauge theory.

\section{Non-simply connected $4$-manifolds}

We focus on just a few questions that the author finds of interest.

\subsection{Particular $4$--manifolds} One of the simplest non-simply connected $4$--manifolds is  $S^{1}\times S^{3}$.
\begin{question}
What is the  classification of group actions on $S^{1}\times S^{3}$?
\end{question}
Jahren and Kwasik \cite{JahrenKwasik2011} have given a thorough study of free involutions on  $S^{1}\times S^{3}$. {\color{black}There are four possible homotopy types for the quotient manifolds. In what is perhaps the most interesting case the quotient is homotopy equivalent to $\mathbb{R}P^{4}\# \mathbb{R}P^{4}$, and there are infinitely many quotients up to homeomorphism. See also Brookman, Davis, and Khan \cite{BrookmanDavisKhan2007} for the latter case and its higher-dimensional analog. The case of free actions of $C_{p}$, $p$ an odd prime (or just odd and square free) has its own special subtleties and has recently been studied by Khan \cite{KhanArXiv2014}.  In particular Khan showed the existence of non-standard actions.

I am unaware of any work on non-free actions on $S^{1}\times S^{3}$.
In this (orientation preserving) case the fixed point set should be a $2$--torus. 
}

It makes sense to study the  rich family of $4$--manifolds that can be expressed as bundles. In particular we suggest study surface bundles.

\begin{problem}
Study group actions on surface bundles over surfaces.
\end{problem}
If both the base and fiber surface have non-positive euler characteristic, then these manifolds are aspherical, having trivial higher homotopy groups. This leads to the next set of considerations.

\subsection{Asymmetric $4$-manifolds}
Aspherical manifolds in all dimensions  figured prominently in the work of Conner-Raymond in their search for manifolds with no compact groups of symmetries. 

\begin{question}
Can a manifold have no nontrivial finite  groups of symmetries?
\end{question}

The first examples of such asymmetric manifolds, due to Bloomberg \cite{Bloomberg1975}, were connected sums of aspherical $4$-manifolds that Conner and Raymond \cite{ConnerRaymond1972} had shown  had the property  that their only  finite  group actions were free. Subsequently asymmetric, aspherical $n$-manifolds  were produced in a few higher dimensions ($7, 11, 16, 22, 29,37$) by Conner, Raymond, and Weinberger \cite{ConnerRaymondWeinberger1972}, and in dimension $3$ by Raymond and Tollefsen \cite{RaymondTollefson1976,RaymondTollefson1982}.

Schultz \cite{Schultz1981} produced infinitely many $n$-manifolds, for any $n\ge 3$, that have no nontrivial  finite  groups of  symmetries, among the class of hypertoral manifolds (admitting a degree one map the the $n$-torus). Initially at least they appear to be connected sums. 

{\color{black}
Subsequently infinitely many asymmetric aspherical $n$-manifolds have been produced in all dimensions $\ge 3$, via hyperbolic geometry. Belolipetsky and Lubotzky \cite{BelolipetskyLubotzky2005} gave a general construction in all dimensions $\ge 2$ of hyperbolic $n$-manifolds with any finite group as isometry group, including the trivial group.  (The cases of dimensions 2 and 3 and of orientation-preserving isometries had been handled earlier and/or independently by others.) Since the center of the fundamental group $\pi$ of a hyperbolic manifold is trivial, it then follows from a result of Borel \cite{Borel1983} (pp. 57--60), that any finite group acting effectively would inject into the (trivial) outer automorphism group $\text{Out}(\pi)$. Thus one needs to produce aspherical manifolds such that $\text{Out}(\pi)$ is trivial.
See Farb and Weinberger \cite{FarbWeinberger2005} for perspectives on Borel's theorem in this context. In the earlier work of Raymond, \emph{et al.}, they proceeded in somewhat the same way, but without input from hyperbolic geometry, and produced instead aspherical fibered manifolds with fundamental group having trivial center and torsion-free outer automorphism group. 
}

As a counterpoint we propose the following:
\begin{problem}
Find a smooth, \underline{simply connected,} closed $4$-manifold with no nontrivial  smooth symmetries.
\end{problem}

This, of course, could be a project for gauge theory. Work  of Edmonds \cite{Edmonds1987} shows that there are no such examples  for topological $4$-manifolds: Every simply connected $4$-manifold  admits  nontrivial topological, locally linear periodic  maps.  In higher dimensions it remains an important problem, pursued most notably  by V.~ Puppe \cite{Puppe2007}.

\subsection{The equivariant Borel conjecture}
Aspherical  manifolds  also play a prominent  role in higher dimensional  topological manifold  theory under the framework of  the Borel and Novikov  conjectures. The Borel conjecture states that a homotopy equivalence between two closed aspherical manifolds  is homotopic to a homeomorphism. The Borel  conjecture has an especially  inviting equivariant version, which we state informally, following Weinberger \cite{Weinberger1994}.

\begin{question}[Na\"ive Equivariant Borel Conjecture]
Is an equivariant homotopy equivalence of aspherical, closed,  topological $G$-manifolds  $G$-homotopic to an equivariant homeomorphism?
\end{question}

This na\"ive version of the equivariant Borel conjecture fails in higher dimensions for relatively simple reasons of not being  properly stratified, and not providing enough normal information to the  fixed point sets. But there are also more subtle sorts of reasons, too. See Weinberger \cite{Weinberger1994}.

The na\"ive equivariant Borel conjecture fails in dimension 4 for actions with two-dimensional fixed points, which may be connected-summed with an action on the sphere with a knotted fixed point $2$--sphere.
\begin{problem}
Find a counterexample to the na\"ive equivariant Borel conjecture in dimension four for pseudofree actions.
\end{problem}

{\color{black}
Connolly, Davis, and Khan  \cite{ConnollyDavisKhan2014} show that a pseudofree involution on $T^4$ is equivalent to the standard involution given by complex conjugation in each factor. In higher dimensions their work gives a complete, nontrivial, classification of pseudofree involutions on tori.

In a follow-up generalization Connolly, Davis, and Khan \cite{ConnollyDavisKhan2015}  give a similar classification of pseudofree actions on aspherical $n$-manifolds, $n\ge 4$, whenever there is an appropriate geometric action to use as a model. Many such model actions are constructed by the process of hyperbolization.

}

As is well-known, every finitely presented group $\pi$ is the fundamental group of a closed orientable $4$-manifold. 

\begin{problem}
If $G$ is a finite group and $\pi$ is a finitely presented group, show that there is a closed, orientable $4$-manifold $M$ with fundamental group $\pi$ on which $G$ acts effectively.
\end{problem}
Perhaps one could even specify the action of $G$ on $\pi$ by a given homomorphism $G\to \text{Out}\, \pi$.

\section{Gauge-theoretic applications to group actions}\label{sec:gauge}
Although the main focus of these notes is topological manifolds, we collect together here some of the main results achieved using gauge-theoretic techniques, including those already briefly mentioned in earlier sections.

In general here the focus is on proving that certain kinds of smooth actions do not exist and, especially, that certain locally linear actions on a smooth manifold are not equivalent to smooth actions. One word of clarification is appropriate here. Many--perhaps all--smooth manifolds admit more than one smooth structure. It is one thing to assert that a certain smooth manifold (such as the standard sphere $S^{4}$ or a standard $K3$ surface) does not admit a smooth action of a certain type. It is something different to say that there is no smooth structure at all on that manifold for which a certain type of action is smooth.

\subsection{Yang-Mills  moduli spaces}
Applications are based on attempting to do as much as possible of Donaldson's construction of an instanton number one moduli space of self-dual (or anti-self-dual) connections equivariantly, producing a $5$-manifold with ideal boundary the given $4$-manifold. One key step causes trouble: applying equivariant transversality to produce an equivariant perturbed, smooth moduli space.  In general equivariant transversality meets obstructions. An  explicit, example of Fintushel (described in Hambleton and Lee  \cite{HambletonLee1992}) shows that nontrivial obstructions can arise in this context. 

Nonetheless there were successes, in which one concentrated on the invariant ASD connections, and only a little bit of normal data. The primary success is with smooth group actions on $4$-manifolds with definite (or trivial) intersection pairing.

Furuta  \cite{Furuta1989}  showed that there are no smooth one-fixed-point actions on (homology) $4$-spheres. Roughly speaking one shows that generically the space of instanton number one $G$-invariant connections forms  a $1$--manifold with an end compactification whose boundary can be identified with the fixed points of the original group action. But a compact $1$--manifold cannot have just one boundary point.

Braam and Matic  \cite{BraamMatic1993}  showed that smooth actions on (homology) $4$-spheres have the same representations at all fixed points. Although they had claimed the result much earlier, it took some years to produce an accepted proof, because of the difficulties with equivariant transversality mentioned above. And in the meantime Hambleton and Lee \cite{HambletonLee1992}  gave a proof using the machinery they developed,. One can assume there are two isolated fixed points. In this case they were able again to study the space of instanton number one $G$-invariant connections and carry along just enough normal data to see that the normal representations must be equivalent.

Both of these were soon followed by more traditional, non-gauge-theoretic proofs, as discussed earlier in this paper.

Other work on (anti-) self-dual moduli spaces and group actions was carried out by Wang \cite{Wang1993} and  Cho  \cite{Cho1990,Cho1991}, who studied the moduli space of invariant connections for manifolds with smooth involution, respectively, with $C_{2}^{k}$ actions.

Hambleton and Lee  \cite{HambletonLee1992} developed an alternative notion of transversality (``equivariant polynomial general position'') via which they were also able to give a gauge-theoretic proof of the result of Edmonds and Ewing \cite{EdmondsEwing1989} on the fixed point data of semifree, pseudofree actions on $\mathbb{C}P^{2}$ and some related matters, also accessible by tradional methods. They also extended the known results beyond what had been accomplished using less powerful techniqes.

More recently, Hambleton and Tanase  \cite{HambletonTanase2004} applied the approach of Hambleton and Lee  \cite{HambletonLee1992}  to prove certain results about smooth actions that are actually false for locally linear actions.  They prove that an action on a connected sum of copies of $\mathbb{C}P^{2}$ has the fixed point data of an equivariant connected sum of actions, which examples of Edmonds and Ewing   \cite{EdmondsEwing1992}  show need not be true in the locally linear case.  They also show that a smooth action on such a manifold that is pseudofree must be semifree. Again an example of Edmonds \cite{Edmonds2005} shows that this is false in general in the locally linear case.

\subsection{Seiberg-Witten theory}
In the 1990's much attention shifted from the Yang-Mills theory to the Seiberg-Witten theory where the moduli spaces in question are generally compact, thus avoiding the thorniest issues associated with Donaldson's theory.  As for group actions, however, results have been slow. The fundamental difficulty of a lack of equivariant transversality is still there. In this context the primary successes have been for smooth group actions on spin $4$-manifolds.

Wang  \cite{Wang1995} gave an early application of Seiberg-Witten invariants when he showed that the quotient of a K3 surface by a free antiholomorphic involution has vanishing Seiberg-Witten invariants.

A major success of Seiberg-Witten invariants was Furuta's famous inequality \cite{Furuta2001}
\[
b_{2}^{+}(X)
 \ge 1- \frac{\sigma(X)}{8}\]
for any smooth, closed, spin $4$-manifold $X$ oriented so that $\sigma(X)\le 0$, with $b_{1}(X)=0$ and $b_{2}^{+}(X)>0$.

While Furuta's theorem is not itself a result about group actions, a key point in the proof is the existence of a certain involution on the moduli space of Seiberg-Witten monopoles.

One of the earliest applications of Seiberg-Witten invariants to more 
general smooth group actions came from Bryan  \cite{Bryan1998}, in which, among other things, he found restrictions on homology representations of smooth involutions (and, more generally, groups of the form $C_{2^{k}}$) on a K3 surface. The approach to the use of the Seiberg-Witten equations was modeled in part on the result of Furuta , and in fact gives an equivariant analog of Furuta's inequality.

\begin{thm}[Bryan  \cite{Bryan1998}] 
Let $X$ be a smooth $4$-manifold with $b_{1}=0$ oriented so that $\sigma(X)\le 0$. Assume $X$ is spin and admits a smoothaction of $C_{2^{p}}$ preserving the spin structure and of odd type. Then under suitable nondegeneracy assumptions
$$b_2^+(X)\geq p+1-\frac{\sigma(X)}{8}$$
\end{thm}

A simple application is the following.

\begin{thm}[Bryan  \cite{Bryan1998}] 
If $C_{2}$ acts smoothly on a K3 surface with isolated fixed points, then $b_{2}^{+}(K3/C_{2})=3$. (In particular, $C_{2}$ must act trivially on the positive definite part of $H_{2}(K3;\mathbb{Q})$.)
\end{thm}

About the same time, Ue \cite{Ue1998} used Seiberg-Witten theory to construct infinitely many smooth free actions on a closed smooth 4-manifold such that their orbit spaces are pairwise homeomorphic, but not diffeomorphic.

\begin{thm}[Ue \cite{Ue1998}] 
For any nontrivial finite group $G$ there exists a smooth 4-manifold that has infinitely many free $G$-actions with the property that their orbit spaces are homeomorphic but mutually nondiffeomorphic.
\end{thm}

In another direction, Kiyono \cite{Kiyono2008} used these ideas to construct nonsmoothable (with respect to any smooth structure), locally linear, pseudofree actions on suitable simply connected spin $4$-manifolds. The current version of his result applies to all but $S^{4}$ and $S^{2}\times S^{2}$. The existence of the actions is an application of work of Edmonds and Ewing \cite{EdmondsEwing1992}.

\begin{thm}[Kiyono \cite{Kiyono2008}]
Let $X$ be a closed, simply connected, spin topological $4$-manifold
not homeomorphic to either $S^{4}$ or $S^{2} \times S^{2}$. Then, for any sufficiently large prime number $p$, there exists a homologically trivial, pseudofree, locally linear action of $C_{p}$ on $X$ which is nonsmoothable with respect to any smooth structure on $X$.
\end{thm}

\begin{problem}
Find an exotic locally linear, pseudofree action on $S^{2}\times S^{2}$. One should seek both actions not equivalent to standard actions and also non-smoothable actions.
\end{problem}

Fang \cite{Fang1998,Fang2001} gave very similar results to Bryan's independently under slightly weaker hypotheses. He also found the following mod $p$ vanishing result for Seiberg-Witten invariants of manifolds with smooth $C_{p}$ actions.

\begin{thm}[Fang \cite{Fang1998}]
Suppose $X$ is a closed, smooth $4$-manifold with $b_1(X) =0$, $b_2^+ (X) >1$, and that $C_p$ acts trivially on $H^{2,+} (X, \bR)$. Then for any $C_p$-equivariant Spin${}^c$-structure $L$ on $X$, the Seiberg-Witten invariant satisfies ${\rm SW}(L) = 0\bmod p$, provided $k_i \le \frac12(b_2^+ -1)$ for $i=0,1, \cdots ,p-1$.
\end{thm}
In the result above the integers $k_{i}=m_{i}-n_{i}$, where $m_{i}$ and $n_{i}$ denote the dimensions of the $\omega^{i}$ eigenspaces of the linear $C_{p}$ actions on $\ker D_{A}$ and $\text{coker }D_{A}$, where $D_{A}$  denotes the Dirac operator corresponding to $L$ and $\omega=\exp(2\pi i/p)$.

J.~H. Kim \cite{Kim2000} also generalized Bryan's result. In particular he developed the analog of Bryan's theorem for actions of $C_{2^{p}}$ of \emph{even} type, with the same inequality as a conclusion.

Both Lee and Li \cite{LeeLi2001} and Bohr \cite{Bohr2002} extended the Bryan results to actions of $C_{2^{p}}$ on non-spin $4$-manifolds with even intersection pairing. The key point is to show that one can lift the action to an action on a covering space that is spin.

In a series of papers Cho  \cite{Cho1999a,Cho1999b,Cho2007} and Cho-Hong 
\cite{ChoHong2002,
ChoHong2003,ChoHong2007}
also investigated conditions imposed on the Seiberg-Witten invariants by a smooth finite group action.

{
\color{black}
The inequality of Hamilton, mentioned at the end of Section \ref{subsec:existence}, uses nonzero Seiberg-Witten classes to bound the number of $2$-spheres in the fixed point set of a smooth action. See Hamilton \cite{Hamilton2016} for applications of the inequality.
}

Recent work of Liu and Nakamura \cite{LiuNakamura2007a} studies homologically nontrivial
 $C_{3}$ actions on a $K3$ surface.  Applying work of Edmonds and Ewing \cite{EdmondsEwing1992}, they give a classification of locally linear actions and then apply Fang's mod $p$ vanishing result to show that certain of these actions cannot be realized smoothly.
 
 Chen and Kwasik \cite{ChenKwasik2011} have also investigated the relationship between symplectic symmetries and the smooth structure on a given symplectice homotopy K3 surface, showing that in a certain sense, an effective action by a large group forces the symplectic homotopy K3 surface to be minimally exotic. Consult the original paper for more statements.
 
Recently, Fintushel, Stern, and Sunukjian \cite{FintushelSternSunukjian2009} have produced infinite families of exotic cyclic group actions on many simply connected smooth $4$-manifolds, all of which are equivariantly homeomorphic but not equivariantly diffeomorphic. The $4$-manifolds in question are constructed as branched cyclic covers of certain simply connected $4$-manifolds $Y$, branched along a smoothly embedded surface $\Sigma$ with the property that the pair $(Y,\Sigma)$ has a nontrivial relative Seiberg-Witten invariant. The branched cyclic cover $X$ has nontrivial Seiberg-Witten invariants, and these are the first such families of manifolds with nontrivial SW invariants. The examples are produced using a variant of the first two authors' ``rim surgery''. These examples are also the first that effectively deal with $2$-dimensional fixed point sets while producing exotic actions, aside from local knotting phenomena arising from counterexamples to the $4$-dimensional Smith Conjecture.

In a similar vein, Kim and Ruberman \cite{KimRuberman2011} have shown how to construct some exotic actions of the cyclic group $C_{mn}=C_{m}\times C_{n}$, where $(m,n)=1$. The actions have the property that the fixed point sets of the two factors are both surfaces intersecting transversely at the isolated global fixed points. The actions are produced by rim surgery starting from standard actions. Topologically the new actions are equivalent to the originals. But smoothly they are distinct, even though the actions of the factors $C_{m}$ and $C_{n}$ are smoothly equivalent to the originals.

{\color{black}
It is an interesting question to study possible bounds on the order of a finite group acting smoothly on a $4$-manifold that does not admit a smooth circle action. Chen \cite{Chen2011} has found such bounds for symplectic or holomorphic actions under suitable technical assumptions.

Recently, however, Chen \cite{Chen2014} has used the knot surgery of Fintushel and Stern  to construct
examples to show that there can be no such bound in the smooth case.}

\appendix

\section{Update of the Kirby Problem List}
Here we briefly record the handful of problems related to group actions on $4$-manifolds that are to be found in the Kirby Problem List, originating at the 1976 Stanford Topology Conference \cite{Kirby1978} and previously updated on the occasion of the 1993 Georgia International Topology Conference \cite{Kirby1993}.

\begin{problem}[Kirby list 4.55]\label{kirby:4.55}
Describe the Fintushel-Stern involution on $S^{4}$ in equations. \textbf{Update (1995):}  No progress.
\end{problem}
\textbf{Update:} No further progress.

\begin{problem}[Kirby List 4.56, P. Melvin]\label{kirby:4.56}
Let $M$ be a smooth closed orientable $4$-manifold  which supports an effective action of a compact connected Lie group $G$. Suppose that $\pi_{1}M$ is a free group.  \textbf{Question:}  Is $M$ diffeomorphic to a connected sum of copies of $S^{1}\times S^{3}$, $S^{2}\times S^{2}$, and $S^{2}{\otimes}S^{2}$? \textbf{Remarks:} The answer to both questions is yes for $G\ne S^{1}\text{ or } T^{2}$; also for $G=T^{2}$ provided the orbit space  of the action (a compact orientable surface) is not a disc with $\ge 2$ holes.  \textbf{Update (1995):}  Still open.
\end{problem}

\textbf{Update:}  No further progress.
\begin{problem}[Kirby list 4.124, Edmonds]\label{kirby:4.124}
(A) Does the fake $\mathbb{CP}^{2}$ (homotopy  equivalent but not homeomorphic to $\mathbb{CP}^{2}$) admit a topological involution?\newline (B) Does the K3 surface admit a periodic diffeomorphism acting trivially on homology?\newline
\textbf{Remarks:} the answer to (B) is no for period 2 [Ruberman, 1995. Matumoto, 1992] . The question is open for odd period.
\end{problem}
\textbf{Update:}\label{prob:K3update}  (A) No known progress. (B) Jin Hong Kim \cite{Kim2008} has used Seiberg-Witten invariants  to show that a smooth homotopy K3 surface does not admit a smooth, homologically trivial, action of $C_{3}$. Experts have expressed doubts, however, about parts of the argument. An alternative proof would be prudent and welcome.
\section{Update of 
problems from the 1984 Transformation Groups Conference}
Here we briefly comment on the status of problems related to $4$--manifolds from the 1984 Boulder conference, edited by Schultz \cite{Schultz1985}.

\begin{problem}[6.10]
What algebraic properties must a knotted $2$--sphere $K$ in $S^{4}$ satisfy if it is fixed (or invariant) under a $C_{m}$ action?
\end{problem}
\textbf{Update:}  No particular progress known.
\begin{problem}[6.11]
What is the complete topological classification of $4$--manifolds with effective $T^{2}$ actions?
\end{problem}
\textbf{Update:}  No further progress beyond the work of Pao \cite{Pao1977a,Pao1977b}. It is the cases with nonsimply connected orbit space that are of interest.
\begin{problem}[6.12]
If $S^{1}$ acts smoothly on a homotopy $4$--sphere $\Sigma$, is $\Sigma$ diffeomorphic to $S^{4}$?
\end{problem}
\textbf{Update:}  The answer is yes by the early work of Fintushel \cite{Fintushel1976} and the positive resolution of the classical Poincar\'e Conjecture by Perelman.
\begin{problem}[6.13]
Does there exist an exotic $\mathbb{R}^{4}$ with no smooth effective finite groups actions?
\end{problem}
\textbf{Update:}  Taylor \cite{Taylor1999} has constructed examples of fake  $\mathbb{R}^{4}$s on which the circle group cannot act smoothly. Taylor also notes that at least some fake $\mathbb{R}^{4}$s do have nontrivial smooth symmetries.
\begin{problem}[6.14]
Are there exotic free $C_{p}$ actions on $S^{1}\times S^{3}$, where $p$ is an odd prime?
\end{problem}
\textbf{Update:}  {\color{black}Khan \cite{KhanArXiv2014} has developed a surgery-theoretic classification of actions on $S^{1}\times S^{n}$ whose existence portion, specialized to  dimension $4$, shows that such exotic actions do exist.}
\begin{problem}[6.15 (Cappell)]
The Cappell-Shaneson example of an exotic smooth $\mathbb{R}P^{4}$ defines a differentiably nonlinear involution on some homotopy $4$--sphere $\Sigma$. Is $\Sigma^{4}$ an exotic homotopy $4$--sphere?
\end{problem}
\textbf{Update:}  
Gompf  \cite{Gompf1991} showed that the simplest one of the Cappell-Shaneson examples does yield an exotic involution on the standard smooth $4$-sphere.  More recently Akbulut  \cite{Akbulut2010} has shown that the members of an infinite family of the Cappell-Shaneson spheres are all diffeomorphic to the standard sphere by showing that they are diffeomorphic to the simplest one, previously dealt with by Gompf. {\color{black}Subsequently Gompf \cite{Gompf2010} has shown that a yet larger family of such examples also yields standard $4$-spheres.}

\nocite{AitchisonRubinstein1984,Akbulut2010,AkbulutKirby1979,AlldayPuppe1993,AnvariHambleton2014,ArtebaniSarti2008,Assadi1990,Baldridge2001,Baldridge2003,Baldridge2004,Bauer1988,BauerWilczynski1996,BelolipetskyLubotzky2005,Bloomberg1975,Bohr2002,Borel1983,BraamMatic1993,Bredon1972,BrookmanDavisKhan2007,Bryan1998,BuchdahlKwasikSchultz1990,CappellShaneson1976,CappellShaneson1979,Carlsson1990,Chen1997,Chen2010,Chen2011,Chen2014,ChenArXiv2011,ChenArXiv2013,ChenArXiv2013a,ChenKwasik2007a,ChenKwasik2008,ChenKwasik2011,ChenKwasikSchultzArXiv2014,Cho1990,Cho1991,Cho1999a,Cho1999b,Cho2007,ChoHong2002,ChoHong2003,ChoHong2007,ConnerRaymond1972,ConnerRaymondWeinberger1972,ConnollyDavisKhan2014,ConnollyDavisKhan2015,DeMichelis1989,DeMichelis1991,deRham1965,Edmonds1981,Edmonds1985,Edmonds1987,Edmonds1988,Edmonds1989,Edmonds1997,Edmonds2005,Edmonds2010,EdmondsEwing1989,EdmondsEwing1992,Fang1998,Fang2001,FarbWeinberger2005,Fintushel1976,Fintushel1977,Fintushel1978,FintushelLawson1986,FintushelPao1977,FintushelStern1981,FintushelStern1984,FintushelStern1985,FintushelSternSunukjian2009,Freedman1984,FreedmanQuinn1990,FukumotoFuruta2000,Furuta1989,Furuta2001,GarbagnatiSarti2007,Giffen1966,Godinho2007,Gompf1991,Gompf2010,Gordon1974,Gordon1986,GuazziMecchiaZimmermann2011,Habegger1982,Hambleton2011,Hambleton2015,HambletonHausmann2011,HambletonKreck1988,HambletonKreck1990,HambletonKreck1993a,HambletonKreck1993b,HambletonKreck1993c,HambletonKreckTeichner1994,HambletonLee1988,HambletonLee1992,HambletonLee1995,HambletonLeeMadsen1989,HambletonMadsen1986,HambletonTanase2004,Hamilton2016,HaynesKwasikMastSchultz2002,Herrera2006,Hsiang1964,Huck1995,HuckPuppe1998,JahrenKwasik2006,JahrenKwasik2011,Jastrzebowski1995,KhanArXiv2014,Kim1993,Kim2000,Kim2007,Kim2008,Kim2011,KimRuberman2011,Kirby1978,Kirby1993,KirbyTaylor2001,Kiyono2008,KiyonoLiu2006,Kottschick2006,Krushkal2005,Kulkarni1982,Kwasik1986a,Kwasik1986b,KwasikLawson1993,KwasikSchultz1988,KwasikSchultz1989,KwasikSchultz1990,KwasikSchultz1991,KwasikSchultz1992,KwasikSchultz1994,KwasikSchultz1995,KwasikSchultz1996,KwasikSchultz1997,KwasikSchultz1999,KwasikVogel1986a,KwasikVogel1986b,Lawson1976,Lawson1993,Lawson1994,LeeLi2001,LiLiu2006,LiLiu2007,LiLiu2008,Liu2000,Liu2005,LiuNakamura2007a,LiuNakamura2008,LongReid2005,Matumoto1986,Matumoto1992,Mazur1961,Mazur1962,Mazur1964,McCooey2002a,McCooey2002b,McCooey2007a,McCooey2007b,McCooey2007c,MecchiaZimmermann2006,MecchiaZimmermann2009,MecchiaZimmermann2011,Melvin1981,Melvin1982,MelvinParker1986,Morimoto1988,Mukai1988,MundetiRieraArXiv2015,Nagami2000,Nagami2001,Nakamura2002,Nakamura2006,Nakamura2009,Nakamura2014,Nikulin1979,Oh1983,Ono1993,Orlik1973,Orlik1982,OrlikRaymond1970,OrlikRaymond1974,Pao1977a,Pao1977b,Pao1978,Park2005,Plotnick1982,Poenaru1960,Puppe2007,RaymondTollefson1976,RaymondTollefson1982,Richardson1961,RuanWang2000,Ruberman1984,Ruberman1995,Ruberman1996,Schultz1980,Schultz1981,Schultz1985,Smith1960,Sumners1975,Sung2015,SungarXiv2011,Szymik2012,Taylor1999,Ue1994,Ue1996,Ue1998,VinogradovKuselman1972,Wang1993,Wang1995,Weinberger1994,Weintraub1975,Weintraub1977,Weintraub1989,Wilczynski1987,Wilczynski1988,Wilczynski1990,Wilczynski1991,Wilczynski1994,WuLiu2014,Yoshida1978,Zeeman1965}

\bibliographystyle{amsplain}
\bibliography{G4Mfds2016}
\end{document}